\documentclass[12pt]{amsart}[1996/10/24]
\relax
\usepackage{amsmath,amscd,stmaryrd,mathrsfs,amssymb}

\advance \topmargin by -\headheight
\advance \topmargin by -\headsep
\evensidemargin \oddsidemargin
\newtheorem{theorem}[equation]{Theorem}
\newtheorem{lemma}[equation]{Lemma}

\newtheorem{proposition}[equation]{Proposition}
\numberwithin{equation}{section}

\newcommand{\Vol}{\operatorname{Vol}}

\newcommand{\Diam}{\operatorname{Diam}}

\begin{document}

\title{A restriction for  singularities on collapsing orbifolds}

\author{Yu Ding}

\address{Department of Mathematics and Statistics,  
California State University, Long Beach, 90840}
\email{yding@csulb.edu}

\date{\today}

\begin{abstract}
Every point $p$ in an orbifold $X$ has a neighborhood that is 
homeomorphic to $G_p\backslash B_r(0)$, where 
$G_p$ is a finite group acting on $B_r(0)\subset{\mathbb R}^n$, so that 
$G_p(0)=0$. 
Assume $X$ is a Riemannian orbifold with isolated singularities that 
is collapsing, that is, $X$ admits a sequence of metrics $g_i$ with 
uniformly bounded curvature, so that for any $x\in X$, 
the volume of $B_1(x)$, with respect to the metric $g_i$, goes to $0$ as 
$i\rightarrow \infty$. 
For such $X$, we prove that $|G_p|\leq (2\pi/0.47)^{n(n-1)}$ for all singularities $p\in X$. 
\end{abstract}

\maketitle

\section{Introduction}
An $n$ dimensional Riemannian 
{\em orbifold}, $X$, is a metric space so that the following is true: 
for any $x\in X$, there exists $r=r(x)>0$ and a 
Riemannian metric $\tilde g_x$ on $B_{2r}(0)\subset{\mathbb R}^n$, a finite 
group $G_x$ (the isotropy group) acting on $(B_r(0), \tilde g_x)$ by isometries, 
so that $G_x(0)=0$, and there is an isometry $\iota_x: B_r(x)\rightarrow 
G_x\backslash B_r(0)$ with $\iota_x(x)=0$. See \cite{Th}. 
We call $x\in X$ a {\em regular point} if $|G_x|=1$; otherwise, 
$x$ is a {\em singular point}. 
We say the curvature of $X$ satisfies 
\begin{equation}
|K_X|\leq \kappa^2,
\end{equation}
if the sectional curvature $K$ of every $(B_r(0), \tilde g_x)$ above 
satisfies $|K|\leq \kappa^2$. We say $X$ 
is {\em collapsing}, if $X$ admits a sequence of metrics, $g_i$, with 
uniformly bounded curvature, so that for any 
$x\in X$, 
\begin{equation}
\lim_{i\rightarrow\infty}\Vol_{g_i}(B_1(x))=0.
\end{equation}

As an example, consider the standard 
${\mathbb Z}_m={\mathbb Z}/m{\mathbb Z}$ action on the sphere $S^2$:

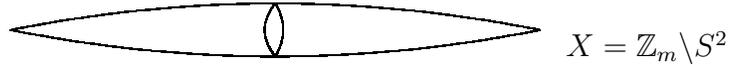
\begin{picture}(400,40)(0,0)
\qbezier(100,20)(200,40)(300, 20)
\qbezier(100,20)(200,0)(300, 20)
\qbezier(200,10)(206,20)(200, 30)
\qbezier(200,10)(192,20)(200, 30)
\put(310,10){$X={\mathbb Z}_m\backslash S^2$}
\end{picture}

\noindent
The quotient orbifold $X={\mathbb Z}_m\backslash S^2$ will be arbitrarily 
collapsed when $m\rightarrow\infty$. However, for any {\em fixed} $m$, 
(i.e. each of the two singularities there is a neighborhood that is 
isometric to ${\mathbb Z}_m\backslash{\mathbb R}^2$, 
where ${\mathbb R}^2$ is equipped 
with some ${\mathbb Z}_m$ invariant metric). 
$X={\mathbb Z}_m\backslash S^2$ can be collapsed only to a certain degree, 
it does not support a sequence of collapsing metrics; in fact, 
otherwise by pull back 
we get a sequence of collapsing metrics on $S^2$, a contradiction to 
the Gauss-Bonnet theorem. 

On the other hand, consider the double of a $2$-dimensional rectangle. 
Clearly it admits a flat metric, thus we obtain a sequence of collapsing 
metrics by rescale.  Notice, in this example, for each of the four 
singularities, the isotropy 
group $G_x$ has order $2$, a quite small number. 

The main result of this paper is 
\begin{theorem}
Assume $X$ is a compact, collapsing orbifold, $p\in X$ is an isolated 
singularity. Then $|G_p|\leq (2\pi/0.47)^{n(n-1)}$.
\label{the-AFO-th-1-GSDG5yJDhjuJDFHj6uDFGHysdrtySFGHy5rSSDH} 
\end{theorem}
If $X$ has an isolated singularity $p$, then the dimension of $X$ 
must be even, and $G_p\subset SO(n)$. 
The bound $|G_p|\leq (2\pi/0.47)^{n(n-1)}$ has its root in the Bieberbach 
theorem of crystallographic groups and Gromov's almost flat 
manifold theorem; however, it fails when we drop the requirement that 
$x$ is an {\em isolated } singularity. For example, we can take {\em any} 
orbifold $X'$ and let 
$X=X\times S^1$; by shrinking the $S^1$ factor, we see $X$ is collapsing 
while there is no restriction on singularities of $X'$. 

Clearly, Theorem 
\ref{the-AFO-th-1-GSDG5yJDhjuJDFHj6uDFGHysdrtySFGHy5rSSDH} 
is a corollary of the following: 
\begin{theorem}
For any $L>0$, there is $\epsilon=\epsilon(n, L)$ so that 
if $X$ is an orbifold with all singularities $q\in X$ satisfying 
$|G_q|<L$, $\Vol(B_1(q))<\epsilon$, then $|G_p|\leq (2\pi/0.47)^{n(n-1)}$ for any 
isolated singularity $p\in X$.
\label{the-AFO-th-2-Hty67ijyufgKJsdfgb5y35SDFGH43sdfsSEFGyhxd} 
\end{theorem}

A {\em nilmanifold}, $\Gamma\backslash N$, is the quotient of 
the (left) action of a discrete, uniform subgroup $\Gamma\subset N$,  
on a simply connected nilpotent Lie group $N$. 
Left invariant vector fields (LIVFs) 
can be defined on $\Gamma\backslash N$.  
An affine diffeomorphism of $\Gamma\backslash N$ 
is a diffeomorphism that maps any local LIVF to some local LIVF. 
In general, a right invariant vector field (RIVF) cannot 
be defined globally in $\Gamma\backslash N$, unless this vector field is 
in the center of the Lie algebra of $N$. However, 
the right invariant vector fields, 
not the left invariant ones, are Killing fields of 
left invariant metrics on $\Gamma\backslash N$. 
An {\em infranil orbifold} is the quotient of a nilmanifold by 
the action of a finite group $H$ of affine diffeomorphisms. If the action 
$H$ is free, we get an {\em infranil manifold}. 

In our previous work, \cite{D}, we generalized 
the Cheeger-Fukaya-Gromov nilpotent Killing structure, \cite{CFG}, 
and the Cheeger-Gromov F-structure, \cite{CG1}, \cite{CG2},  
to collapsing orbifolds. In particular, sufficiently 
collapsed $X$ can be decomposed into a union of {\em orbits}. Each 
orbit ${\mathscr O}_p$ is the orbit of the action of 
a sheaf ${\mathfrak n}$ of nilpotent Lie algebras, which comes from 
local RIVFs on a nilmanifold fibration in the frame bundle $FX$. 
Therefore every 
${\mathscr O}_p$ is an infranil orbifold. The proof of Theorem 
\ref{the-AFO-th-2-Hty67ijyufgKJsdfgb5y35SDFGH43sdfsSEFGyhxd} is based 
on the relation between singularities on $X$ and 
singularities within an orbit ${\mathscr O}_p$ in $X$, as well as 
the nilmanifold fibration on $FX$. 

$X$ is called {\em almost flat}, if 
\begin{equation}
\sup|K_X|^{1/2} \cdot \Diam X\leq \delta_n,
\end{equation}
here $\Diam X$ is the diameter of $X$, $\delta_n$ is a small 
constant that depends only on $n$. 
In \cite{Gr}, Gromov proved that an almost flat {\em manifold} 
$M$ has a finite, normal covering space 
$\tilde M=\Gamma\backslash N$ that is 
a nilmanifold. 
Subsequently, Ruh, \cite{Ru}, proved that $M$ is 
diffeomorphic to 
$\Lambda\backslash N$, where $\Lambda\supset \Gamma$ 
is a discrete subgroup in 
the affine transformation group of $N$. In \cite{Gh2}, Ghanaat generalized 
this to an almost flat orbifold $X$, under the assumption that 
$X$ is {\em good} in the sense of Thurston \cite{Th}, i.e. 
$X$ is the {\em global} quotient of a simply connected manifold $M$. 
There are examples of orbifolds that are not good, see \cite{Th}. In fact, 
without  
much effort, one can remove the assumption that $X$ is good: 

\begin{theorem}
If $X$ is an almost flat orbifold, then $X$ is an infranil orbifold. 

Precisely, there is a nilmanifold 
$\tilde X=\Gamma\backslash N$, a finite group $H$ acting on $\tilde X$ 
by affine diffeomorphism, so that $X$ is diffeomorphic to 
$H\backslash \tilde X$. The order of $H$ is bounded by 
$c_n\leq (2\pi/0.47)^{n(n-1)/2}$. 
Moreover, 
there is a sequence of metrics $g_j$ so that $\Diam (X, g_j)\rightarrow 0$. 
\label{theo-ob-main-ggdfASF5343GSDFGfdasf56HSfhy5}
\end{theorem}
The proof of Theorem 
\ref{the-AFO-th-2-Hty67ijyufgKJsdfgb5y35SDFGH43sdfsSEFGyhxd} does not 
depend on the above result, so we only give a sketch of its proof in 
the Appendix. On the other hand, this implies Theorem 
\ref{the-AFO-th-2-Hty67ijyufgKJsdfgb5y35SDFGH43sdfsSEFGyhxd} for 
almost flat orbifolds immediately, even without the assumption that 
the singularities are isolated.


\section{Proof of Theorem 
\ref{the-AFO-th-2-Hty67ijyufgKJsdfgb5y35SDFGH43sdfsSEFGyhxd}}

If $X$ is an infranil orbifold, then it is easy to obtain 
the bound in Theorem 
\ref{the-AFO-th-1-GSDG5yJDhjuJDFHj6uDFGHysdrtySFGHy5rSSDH}. 
Since the proof contains some ideas for the general case, we give 
full details.   
\begin{lemma}
Assume $X$ is an infranil orbifold. Then $|G_x|\leq (2\pi/0.47)^{n(n-1)/2}$.
\label{cor-gromov-53-order-Gsdfg56y7HDJY65uDRGH36uSRTHQYsGHS}
\end{lemma}
\proof
Assume $X=\Lambda\backslash N$, where $N$ is a simply 
connected nilpotent Lie group, $\Lambda$ is a discrete group of 
affine diffeomorphisms on $N$ so that $X=\Lambda\backslash N$ is compact. 
If $N$ is abelian, then $X$ is a flat orbifold, $\Lambda$ is a 
discrete group of isometries on $N={\mathbb R}^n$ that acts properly 
discontinuously. So the conclusion follows from (the proof of) 
Bieberbach's theorem on crystallographic groups. In fact, it is well known 
that the maximal rotational angle of any $\lambda\in\Lambda$ is either 
$0$ or at least $1/2$. Thus the bound comes from a standard packing 
argument; notice $n(n-1)/2=\text{dim }SO(n)$ and the bi-invariant metric 
on $SO(n)$ has positive curvature. 

We prove the general case by induction on dimension of $X$. 
Remember that $\Lambda$ contains 
a normal subgroup $\Gamma$ of finite index, so that $\Gamma$ is a uniform, 
discrete subgroup of $N$. $X$ the the quotient of 
the $\Lambda/\Gamma$ action on the nilmanifold 
$\tilde X=\Gamma\backslash N$. 
Clearly $G_x$ embeds in $\Lambda/\Gamma$, that is, 
$G_x=\{\bar\lambda\in \Lambda/\Gamma\,|\, \bar\lambda \tilde x=\tilde x\}$; 
here we choose a point $\tilde x$ in $\tilde X=\Gamma\backslash N$ that 
projects to $x\in \Lambda\backslash N$.

Let $C$ be the center of $N$, then $C$ is connected, of positive 
dimension. Since any $\lambda\in \Lambda$ is affine diffeomorphism, 
$\lambda$ moves a $C$-coset in $N$ to a $C$-coset. Therefore 
$\Lambda/\Gamma$ acts on the nilmanifold 
$\tilde X^*=(\Gamma/(\Gamma\cap C))\backslash(N/C)$, the 
quotient $X^*$ is an infranil orbifold of lower dimension. Let 
$\pi: \tilde X\rightarrow \tilde X^*$ be the projection, assume 
$\pi(\tilde x)=\tilde x^*$. Thus we have a homomorphism 
\begin{equation}
h: G_x\rightarrow G_{x^*}.
\end{equation}
$\tilde X$ is a torus bundle over $\tilde X^*$, the fiber is 
$T=(\Gamma\cap C)\backslash C$. Assume $\bar\lambda\in \Lambda/\Gamma$ 
is in $\text{Ker } h$, the kernel of $h$, 
then $\bar\lambda$ fixes every $T$ fiber in $\tilde X$. If, in addition, 
$\bar\lambda$ fixes {\em every point} in the $T$ fiber passing through 
$\tilde x$, 
we claim $\bar\lambda$ must be identity. In fact, on $N$ we have 
$\lambda(z)=a\cdot A(z)$, here $a\in N$ and $A$ is a Lie group 
automorphism of $N$; if $\bar\lambda$ fixes every point in {\em one} 
$T$ fiber, then $A$ is identity on the center $C\subset N$. 
This implies that $\bar\lambda$ is a 
{\em translation} on every $T$ fiber. Since $\bar\lambda$ is of finite 
order and fixes every point in {\em one} $T$-fiber, 
$\bar\lambda$ must be identity. Therefore any element 
$\bar\lambda\in\text{Ker }h$ is decided by its restriction 
on the $T$ fiber passing through $x$; so $\text{Ker }h$ is isomorphic to 
a finite group of affine diffeomorphisms on $T$ that fixes $\tilde x\in T$, 
thus $|\text{Ker }h|$ can be bounded by Bieberbach's theorem. 
Since 
\begin{equation}
|G_x|\leq |G_{x^*}|\cdot |\text{Ker } h|,
\end{equation}
the conclusion follows by induction. 
\qed

In \cite{D}, the existence of nilpotent Killing structure of 
Cheeger-Fukaya-Gromov, \cite{CFG}, is generalized to sufficiently 
collapsed orbifolds. 
We briefly review this construction. 

As in the manifold case, one can define the {\em frame bundle} $FX$ 
of an orbifold $X$. If $B_r(x)\subset X$ is isometric to 
$G_x\backslash B_r(0)$, where $G_x$ is a finite group acting on 
$B_r(0)\subset {\mathbb R}^n$, then locally $FX$ is 
$G_x\backslash FB(0, r)$, here $FB(0, r)$ is the orthonormal 
frame bundle over $B_r(0)$, and $G_x$ acts on $FB(0, r)$ by 
differential, i.e. 
$\tau\in G_x$ moves a frame $u$ to $\tau_* u$. Therefore 
$FX$ is a manifold; strictly speaking, $FX$ is not a fiber bundle. 
Let $\pi: FX\rightarrow X$ be the projection.

Moreover, 
there is a natural $SO(n)$ action on $FX$; on the frames over  
regular points, this $SO(n)$ action is the same one as in the manifold case; 
however, at (the frames over) singular points, this action is {\em not free}. 
As in the work of Fukaya, \cite{F}, see also \cite{D}, any 
Gromov-Hausdorff limit $Y$ of a collapsing sequence $FX_i$ is a manifold. 
Following \cite{CFG}, for 
sufficiently collapsing orbifolds, {\em locally} we have 
an $SO(n)$-equivariant fibration 
\begin{equation}
Z\longrightarrow FX\stackrel{f}{\longrightarrow} Y,
\end{equation}
here the fiber $Z$ is a nilmanifold, $Y$ is a smooth manifold with 
controlled geometry. 

As in \cite{CFG}, we can put a {\em canonical} affine structure 
on the $Z$ fibers, i.e. a canonical way to construct a diffeomorphism from 
a fiber $Z$ to the nilmanifold $\Gamma\backslash N$. In particular, 
there is a sheaf 
${\mathfrak n}$, of a nilpotent Lie algebra of vector fields on $FX$. 
Sections of ${\mathfrak n}$ are local right invariant vector fields 
on the nilmanifold  fibers $Z$. By integrating ${\mathfrak n}$, we get a 
{\em local} action of a simply connected nilpotent Lie group, 
${\mathfrak N}$, on $FX$. Therefore we also call a $Z$ 
fiber an {\em orbit}, we can write $Z=\tilde{\mathscr O}$. 

The fibration $f: FX\rightarrow Y$ is $SO(n)$-equivariant, so 
any $Q\in SO(n)$ moves a $Z$ fiber to a (perhaps another) $Z$ fiber by 
affine diffeomorphism. 
Moreover, the $SO(n)$ action on ${\mathfrak n}$ is {\em locally trivial}, 
that is, if $A\in \mathfrak{so}(n)$ is sufficiently small, then 
$e^A\in SO(n)$ moves a section, ${\mathfrak n}(U)$, of ${\mathfrak n}$ on 
any open set $U\subset FX$, to itself 
(over $U\cap Ue^{A}$). See \cite{CFG} Proposition 4.3. 
In particular, 
the sheaf ${\mathfrak n}$ induces a sheaf, which we also denote by 
${\mathfrak n}$, on the orbifold $X$ away from the singular points. 
An orbit $\tilde{\mathscr O}_{\tilde q}=Z$ 
on $X$ projects down to an {\em orbit} ${\mathscr O}_q$ on $X$. 

Assume $\tilde q\in FX$ is any frame over $q\in X$. Let 
\begin{equation}
I(q)=\{Q\in SO(n)\,|\, Z_{\tilde q}Q=Z_{\tilde q}\}
\end{equation} 
be the the isotropy group of an orbit 
$Z_{\tilde q}={\mathscr O}_{\tilde q}\subset FX$. 
We will simply write $I(q)$ by $I$. Let $I_0$ be 
the identity component of $I$. It can be shown that restricted on 
$Z=\tilde{\mathscr O}_{\tilde q}$, the action of $I_0$ is identical to the 
action of a torus, and the Lie algebra of this torus, $I_0$, 
is in the center of ${\mathfrak n}$. See \cite{CFG}, and \cite{D} 
for more details. Consider the nilmanifold 
\begin{equation}
\check{\mathscr O}_{\tilde q}=\tilde {\mathscr O}_{\tilde q}/I_0. 
\end{equation}
Therefore $Z_{\tilde q}=\tilde {\mathscr O}_{\tilde q}$ is 
a torus bundle over $\check{\mathscr O}_{\tilde q}$. 
Notice, on $Z=\tilde{\mathscr O}_{\tilde q}$, $I$ moves $I_0$ fibers to 
$I_0$ fibers, thus  
the orbit ${\mathscr O}_q$ is the quotient of 
$\check{\mathscr O}_{\tilde q}$ by the 
action of the finite group $I/I_0$. Therefore ${\mathscr O}_q$ is 
an infranil orbifold. In particular, the singularities within 
${\mathscr O}_q$ satisfies the bound in Lemma 
\ref{cor-gromov-53-order-Gsdfg56y7HDJY65uDRGH36uSRTHQYsGHS}. 

It is important to remark that the above structure is not trivial: 
\begin{lemma}
Let $L$ be any integer. Then there is $\epsilon=\epsilon(n, L)$, so that 
if $X$ is an orbifold with $|G_x|\leq L$, $\Vol(B_1(x))\leq \epsilon$ for 
all $x\in X$, then every ${\mathfrak n}$-orbit 
${\mathscr O}$ on $X$ is of positive dimension.
\label{lemma-gromov-positive-dim-gsdfg56SDFHG345NHdtnju767uBDSTR24}  
\end{lemma}
\proof (sketch)
For any unit vector  
$A\in \mathfrak{so}(n)$, the bound in $|G_x|$ implies that $e^{tA}$ does not 
have fixed point in $\tilde{\mathscr O}_{\tilde q}$ unless $t=0$ or 
$|t|>cL^{-1}$. However, for sufficiently collapsed orbifolds, there is a 
vector $B$ in the center of ${\mathfrak n}$ so that $B$ generates a closed 
loop in $\tilde{\mathscr O}_{\tilde q}$ 
that is shorter than $cL^{-1}$, 
therefore $B$ cannot be in the Lie algebra of 
$I_0$, which is in both $\mathfrak{so}(n)$ and the center of ${\mathfrak n}$. 
Thus the orbit $\tilde{\mathscr O}_{\tilde q}$ is not contained in a single 
$SO(n)$ orbit in $FX$, 
so ${\mathscr O}_q$ is of positive dimension in $X=FX/SO(n)$. 
See \cite{D} for more details. 
\qed


\proof[Proof of Theorem 
\ref{the-AFO-th-2-Hty67ijyufgKJsdfgb5y35SDFGH43sdfsSEFGyhxd}]

Assume $p\in X$ is an isolated singular point, $\tilde p\in\pi^{-1}(p)$ 
is in $FX$.  
$Z_{\tilde p}=\tilde{\mathscr O}_{\tilde p}$ is the fiber that projects 
to ${\mathscr O}_p$. 
Let $I(p)$, $I_0$, $\check{\mathscr O}_{\tilde p}$ be as above. Let 
\begin{equation}
K_p=\{Q\in SO(n)\,|\, \tilde p Q=\tilde p\}.
\end{equation}
Thus $K_p$ is a subgroup of $I$, and $|K_p|=|G_p|$.  
Let 
\begin{equation}
K_p^-=\{Q\in K_p\,|\,(\tilde p'I_0) Q=(\tilde p'I_0)\,\,\, \text{ for all }
\tilde p' \in Z_{\tilde p}\}.
\end{equation}
Thus $K_p^-$ is a normal subgroup of $K_p$. 

\begin{lemma}
If $Q\in K_p$ fixes every point in $Z_{\tilde p}$, then $Q$ 
is the identity in $SO(n)$.
\end{lemma}
\proof
Potentially $Q$ may fix every point in $Z_{\tilde p}$ while moving some 
points of $FX$ that are outside $Z_{\tilde p}$. We will rule out this 
possibility. 

By assumption, $p$ is an isolated singularity. 
For any $Q$ that is not identity, the  
connected component of the fixed point set of $Q$ that passes through 
$\tilde p$ must project to $p$ under $\pi: FX\rightarrow X$, because 
away from $\pi^{-1}(p)$ the $SO(n)$ action is free. Therefore 
$\pi(Z_{\tilde p})=p$ is a single point in $X$, this 
contradicts the fact that the ${\mathfrak n}$-orbits on $X$ are of 
positive dimension; see Lemma 
\ref{lemma-gromov-positive-dim-gsdfg56SDFHG345NHdtnju767uBDSTR24}.  
\qed

In particular, we have a {\em faithful} representation of $K_p$ in the 
affine group of $Z_{\tilde p}$, i.e. we can identify $K_p$ with 
the {\em restricted action} of the group $K_p$ {\em on} $Z_{\tilde p}$. 

Take any $Q\in K_p^-$ that is not identity in $SO(n)$. By definition 
$Q$ fixes $\tilde p$. 
If $Q$ fixes {\em every} point in $\tilde p I_0$, as in Lemma 
\ref{cor-gromov-53-order-Gsdfg56y7HDJY65uDRGH36uSRTHQYsGHS}, $Q$ is a 
translation on {\em every} $I_0$ fiber;
because $Q$ is of finite order and $Q$ moves every $I_0$ fiber to itself, 
$Q$ necessarily fixes every point in $Z$; thus $Q$ is identity. 
So $Q$ rotates the tangent plane of $\tilde pI_0$ at $\tilde p$. Therefore 
$K_p^-$ is isomorphic to a finite group of affine diffeomorphisms on the torus 
$\tilde p I_0$. 
By the Bieberbach theorem, 
\begin{equation}
|K_p^-|\leq (2\pi/0.47)^{k(k-1)/2}, \ \ \ \ k=\text{dim } I_0. 
\end{equation}
Recall that Bieberbach's theorem implies that all finite subgroups of 
$SL(n, {\mathbb Z})$ have 
a uniform upper bound in order.

We have 
\begin{equation}
{\mathscr O}_p=H\backslash(Z_{\tilde p}/I_0), 
\end{equation}
where $H=I/I_0$ is a finite group. 
Let $H_p$ be the subgroup of $H$ that fixes $p$. 
Now we get an {\em embedding} 
\begin{equation}
K_p/K_p^- \subset H_p.
\end{equation}
By Lemma \ref{cor-gromov-53-order-Gsdfg56y7HDJY65uDRGH36uSRTHQYsGHS}, 
\begin{equation}
|H_p|\leq (2\pi/0.47)^{i(i-1)/2}, \ \ \ \ i=\text{dim } {\mathscr O}_p.
\end{equation}
Thus 
\begin{equation}
|G_p|=|K_p|=|K_p/K_p^-|\cdot |K_p^-|\leq (2\pi/0.47)^{n(n-1)}.
\end{equation}
\qed


\section*{Appendix: orbifold version of Gromov's theorem}

Assume $M$ is an almost flat manifold.
The proof of Gromov's almost flat 
manifold theorem has two natural steps. 
First, by a delicate calculus of 
geodesic loops, one gets a finite, normal covering $\hat M$ of $M$; see 
\cite{Gr}, \cite{BK} Chapters 1-3. Next, one uses the deformation technique 
of Ruh \cite{Ru} and Ghanaat \cite{Gh} to show that $\hat M$ is a 
nilmanifold. We follow closely this approach. See also 
\cite{Ro}.

Rescale $X$ so that $\Diam X=1$, therefore the curvature of $X$ is bounded 
by $\delta_n^2$. 
Let $p\in X$ be a regular point, and $v$ be any unit tangent vector at $p$. 
The geodesic $\gamma$ so that $\gamma(0)=p$ and $\gamma'(0)=v$ is defined 
at least on a short interval $[0, \delta)$. One can extend 
it to be a geodesic $\gamma: [0, \infty)\rightarrow X$. In fact, 
the extension exists as long as $\gamma(t)$ is a regular point;  
if $\gamma(t)=x$ is singular, 
recall there is a neighborhood $U_x$ of $x$, 
an Euclidean ball $V_x$ with metric $\tilde g$ 
so that $U_x=G_x\backslash V_x$; 
we can lift $\gamma|_{(t-\epsilon,\, t]}$ 
to a geodesic on $V_x$, extend it beyond $t$ on $V_x$ 
and then project to $U_x\subset X$. Therefore, we get a map 
\begin{equation}
\exp: {\mathbb R}^n=T_p \longrightarrow X,
\end{equation}
so that $\exp(sv)=\gamma(s)$. This map is called the {\em develop map}; see 
\cite{Th}, \cite{CT}. It plays the role of exponential map. Let 
\begin{equation}
\rho< r_{\text{max}}=\pi/\sqrt{\delta_n}. 
\end{equation}
By standard comparison theorem, $\exp$ is {\em nonsingular} on 
$B_{\rho}(0)$. When $x=\exp(tv)$ is a singular point this means,  
there is a neighborhood $\hat V$ of $tv\in {\mathbb R}^n$ 
so that the following 
local lift $\hat\exp$ is nonsingular (it is always smooth): 
\begin{equation}
\hat V \stackrel{\hat\exp}
\longrightarrow V_x \stackrel{\text{pj}_x}\longrightarrow U_x\subset X, 
\ \ \ \ \ \ \ \ 
\exp|_{\hat V}=\text{pj}_x\circ \hat\exp. 
\end{equation}
In particular, 
the Riemannian metric $g$ on $X$ pulls back to be a smooth 
metric $\tilde g$ on $B_\rho(0)\subset {\mathbb R}^n$. 
In the following, 
we will assume that $B_\rho(0)\subset T_p$, write $T_p$ 
for simplicity, is equipped 
with the metric $\tilde g$.   

We use the convention that for two loops $\alpha, \beta$ in $X$, 
the product $\alpha\cdot\beta$ is the loop 
that goes along $\beta$ first, then $\alpha$. Define
\begin{equation}
\pi_\rho=\{\text{ geodesic loops } \gamma \text{ at } p\in X   
\text{ with length }
|\gamma|\leq\rho \}.
\end{equation}
Choose $\rho$ so that $2\Diam(X)\leq 2 \ll \rho \ll r_{\text{max}}$. 
Assume $\alpha, \beta$ are two loops at $p$, 
$|\alpha|+|\beta|<r_{\max}$. Then we can lift $\alpha\cdot\beta$ to a 
piecewise geodesic on $T_p$; such a path is homotopic (rel. end points) 
to a unique minimal geodesic $\tilde\gamma$ on $T_p$. 
Their {\em Gromov product}, 
$\alpha*\beta$, is defined to be $\exp(\gamma)$. 
This kind of homotopy on $X$, which can be lifted to $T_p$ through 
$\exp$, is called {\em Gromov homotopy}. 
In general $\alpha*\beta$ is a geodesic loop in $X$; 
but {\em not} the minimal 
geodesic in the homotopy class of $\alpha\cdot\beta$ - examples 
are abundant on simply-connected orbifolds, e.g. the double of a 
standard square. 

By lifting to $(T_p, \tilde g)$, we see that 
$\alpha*\beta$ is homotopic to 
$\alpha\cdot\beta$ 
via a {\em length decreasing} homotopy rel. end points. The 
assumption $|\alpha|+|\beta|<r_{\max}$ guarantees that the homotopy 
takes place in a domain of $T_p$ on which $\exp_p$ is nonsingular. 
This is the orbifold version of the {\em short homotopy} in \cite{BK}. 

Put the standard affine structure on ${\mathbb R}^n=T_p$. 
For each loop $\gamma\in\pi_\rho$ at $p$, we define the 
{\em holonomy motion}:
\begin{equation}
m(\gamma):T_p\rightarrow T_p, \ \ 
v\mapsto r(\gamma)v+t(\gamma).
\end{equation}
Here $r(\gamma)\in O(n)$ is the parallel translation along $\gamma$, 
$t(\gamma)=|\gamma|\cdot \gamma'(|\gamma|)$.
In another word, $m(\gamma)v=A(|\gamma|)$, where 
$A$ is the solution of the following ODE along $\gamma$:
\begin{equation}
\left\{\begin{array}{ccc}
\nabla_TA=T, \\
A(0)=v.
\end{array}\right.
\end{equation}
In particular, $m(\alpha\cdot\beta)=m(\alpha)\circ m(\beta)$. 
As in \cite{BK}, one can 
estimate the error caused by homotopy and conclude that 
\begin{equation}
m(\alpha*\beta)=m(\alpha)\circ m(\beta) \ \ \ \ \text{with error up to }\,\,
c\Lambda^2\rho^2 < \epsilon_n\ll 1.  
\label{eq-ob1-homo-error-GGFsdgrgNDghj7udfHSTRHjw4yjysdSAEGRsw}
\end{equation}

The following result of Gromov, \cite{Gr}, is of 
fundamental importance: 

\begin{proposition} 
There is a positive constant $\theta\ll 1$, so that 
if $|\alpha|\leq {\rho/3}$ and $|r(\alpha)|\leq 0.48$, then 
\begin{equation}
\text{maximum rotational angle of }\, r(\alpha)
\leq \theta \rho^{-1}|\alpha|.
\end{equation}
\end{proposition}

In Chapter 3 of \cite{BK}, a detailed proof of the above proposition 
were given in the manifold case. All the arguments there (which 
involves Jacobi field estimates) were carried out on the tangent space 
$T_p$; therefore the proof in \cite{BK} can be used here, by   
replacing the topological Gromov product in \cite{BK} with the 
one we defined above. 

Thus the rotation parts of $m(\gamma)$, with $\gamma$ not too long,  
are concentrated in finitely many regions, $S_1, S_2, ...,S_N \subset O(n)$, 
which are mutually separated by distance at least $0.47$. We can 
put these loops into finitely many classes, loops in the same class 
have rotational parts in the same $S_i$. There is 
a natural group structure on the set $H$ of these classes, and 
$H$ can be embedded into $O(n)$, with  
\begin{equation}
|H|\leq (2\pi/0.47)^{\text{dim} O(n)}.
\end{equation}

Fix an orthonormal frame $u$ at $0\in (T_p, \tilde g)$, 
parallel translate it 
along normal geodesics, we get an almost 
parallel frame field $u$ on $B_{\rho}\subset T_p$. 
Let $q\in X$ be a regular point. Let $\tilde q_1, \tilde q_2\in T_p$ 
so that $\exp(\tilde q_1)=\exp(\tilde q_2)=q$. Then there is an element 
$h\in H$, so that 
\begin{equation}
d(\exp_*(u(\tilde q_1))h, \exp_*(u(\tilde q_2)))\leq \epsilon_n\ll 0.47.
\end{equation} 
In particular, if $\{\tilde q_1, \tilde q_2, ...\}= B_\rho(0)\cap\exp^{-1}(q)$, 
then the frames $\exp(u(\tilde q_i))$, $i=1, 2, ...$,  
also concentrated in regions $S_1^q, ...,S_N^q$. There is a 
right action by the discrete group $H$, which permutes these regions, 
with error up to $\epsilon_n\ll0.47$. 
By a standard center of mass argument, for each $h\in H$ we get a 
frame $C_q(h)$ at $q$ so that 
\begin{equation}
C_q(h)h'=C_q(hh'), \ \ \ \ \text{for any }\,\, h'\in H. 
\end{equation} 
Consider the set $X'\subset FX$ of all frames $C_qh$, denote by $\tilde X$ 
its closure. So $\Diam\tilde X\leq 2|H|$. Now it is easy to see that 
$\tilde X$ is a smooth manifold, and $X=H\backslash \tilde X$, 
see \cite{BK}. Write $\pi:\tilde X\rightarrow X$ the projection.  
Any point $\tilde q\in\tilde X$ can also be viewed as a frame 
$v(\tilde q)$ on $X$. Then 
$\pi^*v(\tilde q)|_{\tilde q}$ 
gives a {\em global} orthonormal frame field on $\tilde X$. 

Finally, following the deformation technique in Ruh, \cite{Ru}, and  
Ghanaat \cite{Gh}, we see $\tilde X$ is a nilmanifold, 
that concludes the proof of Gromov's theorem.

{\sc Acknowledgements}. 
We are grateful to Prof. Tian for very helpful suggestions.

\end{document}